\documentclass{siamart251216}

\usepackage{xurl}
\usepackage{amssymb}
\usepackage{algorithm}
\usepackage{algpseudocode}
\usepackage{booktabs}
\usepackage{graphicx}
\usepackage{bm}
\usepackage{xcolor}
\usepackage[normalem]{ulem}
\usepackage{epstopdf}
\usepackage[section]{placeins}
\usepackage{tikz}
\usepackage{siunitx}

\usepackage{comment}

\graphicspath{{}{Figs/}}

\newcommand{\csrc}{Beijing Computational Science Research Center, Beijing, China}
\newcommand{\odu}{Department of Mathematics \& Statistics,
  Old Dominion University,\\
  Norfolk, VA 23508, USA}

\usepackage{unicode-math-input}

\newcommand{\Rey}{\mathrm{Re}}
\newcommand{\May}{\mathrm{Ma}}

\newcommand{\cs}{c_{\mathrm{s}}}

\newcommand{\ssstyle}{\scriptscriptstyle}
\newcommand{\fazero}{f_\alpha^{\ssstyle (0)}}

\newcommand{\tensorfont}[1]{\bm{\mathsf{#1}}}
\newcommand{\tensorgreek}[1]{\pmb{#1}}
\newcommand{\Adelta}{\tensorfont{A}_{\delta}}
\newcommand{\Tdelta}{\tensorfont{T}_{\!\!\delta}}
\newcommand{\qvector}[1]{\bm{\mathrm{#1}}}
\newcommand{\dtH}{\delta_t^{\mathrm{H}}}
\newcommand{\dth}{\delta_t^{\mathrm{h}}}
\newcommand{\dxH}{\delta_x^{\mathrm{H}}}
\newcommand{\dxh}{\delta_x^{\mathrm{h}}}

\newcommand{\ie}{i.e.\@}

\newcommand{\Reygrid}{\Rey^{*}}
\newcommand{\omegaplus}{\omega_+}
\newcommand{\omegaminus}{\omega_-}
\newcommand{\weightLBE}{\lambda}
\newcommand{\gammaPLBE}{\gamma}
\newcommand{\Uzero}{U_{0}}
\newcommand{\Ulid}{U_{\mathrm{lid}}}
\newcommand{\dtLB}{\delta_t}
\newcommand{\dxLB}{\delta_x}
\newcommand{\viscLB}{\nu}
\newcommand{\Mapre}{\May_{\mathrm{p}}}
\newcommand{\Rmap}{\tensorfont{R}}
\newcommand{\Phistate}{\tensorgreek{\Phi}_{p}}
\newcommand{\Phires}{\tensorgreek{\Phi}_{b}}
\newcommand{\Pcor}{\tensorfont{P}_{b}}

\newcommand{\fplus}{f_{\alpha}^{+}}
\newcommand{\fminus}{f_{\alpha}^{-}}

\newcommand{\faeqpm}{f_{\alpha}^{\pm\!\ssstyle (0)}}
\newcommand{\feqplus}{f_{\alpha}^{+\!{\ssstyle (0)}}}
\newcommand{\feqminus}{f_{\alpha}^{-\!{\ssstyle (0)}}}
\newcommand{\FASdefect}{\qvector{D}}
\newcommand{\smoothpre}{n_{\mathrm{pre}}}
\newcommand{\smoothpost}{n_{\mathrm{post}}}
\newcommand{\smoothcoarse}{n_{\mathrm{c}}}
\newcommand{\qf}{\qvector{f}}
\newcommand{\qu}{\bm{u}}
\newcommand{\qxi}{\bm{\xi}}
\newcommand{\qcc}{\bm{c}}

\newcommand{\vrj}{\bm{r_j}}
\newcommand{\phat}{\hat{p}}
\newcommand{\falpha}{f_\alpha}
\newcommand{\fbaralpha}{f_{\bar\alpha}}
\newcommand{\walpha}{w_\alpha}
\newcommand{\Ncavity}{N}
\newcommand{\clattice}{c}

\title{A Nonlinear Multigrid Method for the Preconditioned Lattice Boltzmann Equation}

\author{%
Xingzhou Ji\thanks{
\csrc\ (\email{jixz@csrc.ac.cn}).
}
\and
Li-Shi Luo\thanks{
\odu\ (\email{lluo@odu.edu}).
Corresponding author.
}
}

\begin{document}
\maketitle

\begin{abstract}
This work developed a preconditioned lattice Boltzmann equation
(PLBE) smoother for the residual-driven full approximation scheme
(FAS) multigrid to accelerate steady-state LBE calculations.
To overcome the deterioration of LBE smoothing at low viscosity,
%Motivated by the relaxation-dependent deterioration of the LBE
%smoothing, 
we used a two-relaxation-time (TRT) collision model together with
characteristic-velocity and kinetic preconditioning to keep the relaxation
frequency for the even-order moments within a favorable range, while
independently selecting the relaxation frequency for the odd-order moments
to reduce their local amplification factor.
We further investigated F-, W-, and hybrid V/W-cycles to address residual stagnation
observed with V-cycles at high Reynolds numbers.
The 2D lid-driven cavity flow was used as a test case.
At $\Rey=5000$, TRT-PLBE-FAS
%with TRT collision model reduces solution time by about 17\% 
%compared with PLBE-FAS with the collision 
%of single relaxation time (SRT); and 
achieves an approximately 27-fold speedup over single-grid TRT-PLBE.
%reaches a common residual threshold 27 times faster than single-grid TRT-PLBE???.
At $\Rey=20000$ on a $2048^2$ grid, the hybrid V/W-cycle
reduces solution time by 14\% compared with the W-cycle.
%The computed solutions agree well with high-resolution
%reference results up to $\Rey=20000$.
\end{abstract}



\section{Introduction}
\label{sec:trt-introduction}

The lattice Boltzmann equation (LBE) is widely used for the simulation
of nearly incompressible flows owing, in part, to its
algorithmic simplicity. The LBE core algorithm consists of
local collision and streaming between neighboring nodes
on a Cartesian grid \cite{lallemand2021lattice}.
%Collision is local to each lattice node, whereas streaming involves
%regular data movement between neighboring nodes on a Cartesian grid.
The LBE algorithm requires relatively few floating-point operations
(FLOPs) per iteration and consequently has a low computational cost;
it also
%These features lead to a low computational cost per iteration and
facilitates efficient implementation on modern parallel architectures.
For steady-state calculations, however, the primitive LBE iteration
may require a large number of steps to eliminate slowly decaying
errors; thus, its inexpensive iterations do not necessarily lead to
%translate into
an efficient way to obtain a steady-state solution.
%a short time to a converged steady solution.

Geometric multigrid methods offer
%offers 
a natural
%natural 
way to
%address 
alleviate this impediment.
%limitation.
The regular Cartesian lattice inherent to the LBE
naturally facilitates the construction of a hierarchy of
coarser grids, while its simple and hence inexpensive collision-streaming
iteration may serve directly as a nonlinear smoother
without requiring a separate local linear solver.
The full approximation scheme (FAS) is particularly suitable
because it acts directly on the nonlinear problem and
transfers corrections between grid levels.
Nevertheless, the low cost of an LBE smoothing step alone
does not ensure
%ensure 
efficient multigrid convergence:
the smoother must damp errors that are not adequately
represented on coarser grids, while the coarse-grid
correction must reduce the remaining error components.

Multigrid acceleration of the LBE has been investigated previously. Mavriplis \cite{mavriplis2006multigrid} developed a FAS multigrid based directly on the collision--streaming LBE, demonstrating nearly grid-independent convergence at a fixed Reynolds number and achieving a significant speedup over the single-grid LBE for two-dimensional lid-driven cavity flow at $\Rey=100$. However, the multigrid convergence rate deteriorated significantly as the Reynolds number increased at fixed grid resolution. Armstrong and Peng \cite{Armstrong2019ccp} extended this approach to the multiple-relaxation-time (MRT) collision model and three-dimensional flows, demonstrating improved numerical stability and substantial speedups over the single-grid  LBE. Their investigations highlighted the convergence benefits of W-cycles and red--black Gauss--Seidel smoothing. Nevertheless, the improved convergence of Gauss--Seidel did not translate into comparable reductions in CPU time for three-dimensional problems, primarily because of poorer memory locality and additional storage requirements. Gsell \emph{et al.} \cite{gsell2020multigrid} developed a dual-time-stepping method for unsteady flows, employing a TRT-based W-cycle FAS multigrid solver to solve the pseudosteady problem arising at each physical time step. However, the potential of tuning TRT relaxation parameters to improve multigrid smoothing was not systematically investigated. Overall, previous studies have largely emphasized more extensive coarse-grid cycling and stronger iterative smoothers to accelerate FAS-LBE convergence. Although effective in certain configurations, W-cycles involve repeated coarse-grid visits, increasing computational work and potentially amplifying communication overhead on massively parallel architectures \cite{AvnatYavneh2023}. Multi-color Gauss--Seidel smoothing allows parallel updates within each color but requires sequential updates across colors and may suffer from reduced memory locality \cite{Armstrong2019ccp,BakerEtAl2011}. Moreover, these approaches have not fully resolved the deterioration of multigrid convergence with increasing grid Reynolds number. These limitations motivate the development of FAS-LBE methods that retain the parallel-friendly collision--streaming structure while exploiting the intrinsic degrees of freedom in the LBE collision and preconditioning models to improve nonlinear smoothing and mitigate Reynolds-number-dependent convergence deterioration.

In a recent work \cite{JiLuo2026cmame}, we formulated the LBE as a
nonlinear residual defined by the complete collision--streaming mapping.
This residual-driven formulation separates the discrete LBE operator
from the nonlinear solution algorithm and provides a common framework for
fixed-point iteration, FAS multigrid, and Newton--Krylov methods.
The FAS implementation primarily employed single-relaxation-time
(SRT) collision, with the preconditioned LBE (PLBE) also considered
as a nonlinear smoother. Although the resulting method accelerated
steady-state calculations, its smoothing behavior remained sensitive
to the collision relaxation frequency. The smoothing properties of 
the SRT-LBE depend strongly on the relaxation frequency and may 
deteriorate when it is either too small or too close to its upper admissible limit. 
Moreover, the use of a single relaxation frequency for all moments 
may lead to insufficient damping of higher-order kinetic modes 
whose relaxation rates are not constrained by the physical viscosity.

The objective of the present work is to develop a strategy to mitigate
or overcome the Reynolds-number-dependent deterioration
of smoothing in the LBE multigrid algorithm \cite{JiLuo2026cmame}.
%The present work develops a two-relaxation-time preconditioned LBE
%(TRT-PLBE) smoother within the residual-driven FAS framework. 
%To address the relaxation-dependent deterioration of smoothing
%identified in our previous work \cite{JiLuo2026cmame}, 
We propose a systematic strategy for selecting relaxation parameters
to improve nonlinear smoothing using the preconditioned LBE
with a two-relaxation-time (TRT) collision model. 
The TRT model relaxes the even- and odd-order moments of the distribution 
functions separately with different relaxation frequencies
%The TRT collision allows the even and odd components of the
%distribution functions to relax independently
\cite{Ginzburg2008cicp427}, 
while kinetic preconditioning introduces an additional parameter
through the equilibrium distribution \cite{guo2004preconditioned}.
%We exploit 
These degrees of freedom are to be exploited 
in two complementary ways:
the characteristic velocity and preconditioning parameter are
adjusted to confine
%place 
the even relaxation frequency in a range
favorable for nonlinear smoothing at a prescribed Reynolds
number, while the odd relaxation frequency 
is selected to suppress or eliminate the local
amplification factor for the odd-order nonconserved moments.
The proposed approach aims to mitigate Reynolds-number-dependent
smoothing deterioration while
%These improvements are achieved while 
retaining the original collision--streaming algorithmic structure of
the discrete LBE within the residual-driven LBE framework
\cite{JiLuo2026cmame}.

%Improved 
It should be noted that
effective smoothing alone does not
%, by itself, 
guarantee rapid convergence of the complete multigrid iteration.
In particular, some high-Reynolds-number V-cycle tests may suffer
%exhibit 
residual stagnation, suggesting that
additional coarse-grid work is needed to further reduce the remaining 
slowly converging error components. Therefore 
%examine more extensive
coarse-grid traversals will be studied more extensively, 
including the PETSc F-cycle, W-cycles, and hybrid V/W-cycles. 
%The F-cycle mode includes a coarse-to-fine startup stage 
%followed by V-cycle iterations; 
%it is not treated as a repeated ??? classical F-cycle. 
We therefore investigate the relative efficiency of these cycle
strategies through numerical tests.

The remainder of this paper is organized as follows.
\Cref{sec:trt-plbe-fas} presents the preconditioned 
LBE with the TRT collision model,
the residual-driven nonlinear smoother, the FAS construction, and
the selection of relaxation parameters across grid levels.
\Cref{sec:trt-numerical-results} presents the numerical results. 
We first compare the different smoothing configurations,
single-grid and multigrid convergence, and multigrid cycle strategies.
The resulting PLBE multigrid algorithm is validated using the
2D lid-driven cavity flow at high Reynolds numbers.

%followed by comparisons of high-Reynolds-number cavity solutions with
%reference data.

\section{PLBE-FAS method}
\label{sec:trt-plbe-fas}

The residual-driven lattice Boltzmann equation (LBE) \cite{JiLuo2026cmame} 
is built upon the complete collision-streaming mapping as its nonlinear 
discrete operator.  We combine the full approximation scheme (FAS) 
with the preconditioned LBE using the TRT collision model on a Cartesian grid.
%two-relaxation-time (TRT) collision and the preconditioned 
%LBE (PLBE) equilibrium on a hierarchy of uniform Cartesian grids.

\subsection{TRT collision and PLBE equilibrium}
\label{subsec:trt-collision}

\begin{figure}[htbp]
\centering
\includegraphics[width=0.25\linewidth]{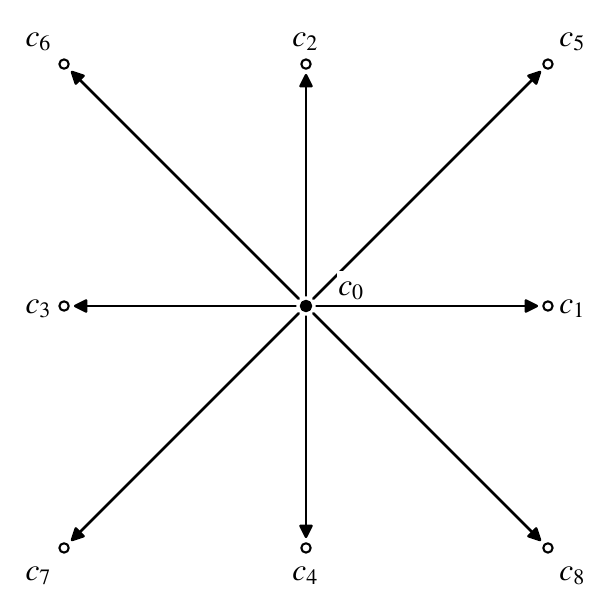}
\caption{Schematic of the D2Q9 discrete velocity set.
The arrows represent the eight nonzero discrete velocities,
while the central dot corresponds to the zero velocity.}
\label{fig:d2q9-velocities}
\end{figure}

We consider a \emph{cell-centered} uniform Cartesian grid $\dxLB \mathbb{Z}_d$
in two dimensions ($d=2$) with the grid spacing $\dxLB$.
Let $\vrj$ denote the position of a cell center indexed by $\bm{j} \in \mathbb{Z}_d$,
and $t_n := n \dtLB$ the discrete time
associated with the LBE evolution.
In this work, we consider the two-dimensional nine-velocity
(D2Q9) LBE model.
%For the D2Q9 model, 
The vector of the local distributions is
$\qf := (f_0,\ldots,f_8)^\dagger$, where
$\falpha := \falpha(\vrj,\, t_n)$ is the distribution associated with the
discrete particle velocity $\qxi_\alpha$, $\alpha=0,\, 1,\ldots,8$.
The discrete velocities are defined by
$\qxi_\alpha := \clattice \qcc_\alpha$, where
$\clattice := \dxLB/\dtLB$ is the lattice speed unit,
$\dtLB$ is the discrete time step size. 
The dimensionless D2Q9 velocity vectors $\qcc_\alpha$ are
\begin{equation}
\begin{aligned}
\qcc_0 := \begin{pmatrix} 0 \\ 0 \end{pmatrix},
\quad
%\(0,0),\\
&\qcc_1 := \begin{pmatrix} 1 \\ 0 \end{pmatrix},
\quad
\qcc_2 := \begin{pmatrix} 0 \\ 1 \end{pmatrix},
\quad
\qcc_5 := \begin{pmatrix} 1 \\ 1 \end{pmatrix},
\quad
\qcc_6 := \begin{pmatrix} -1 \\ 1 \end{pmatrix},
\\
&\qcc_3 := -\qcc_1,
\quad
\qcc_4 := -\qcc_2 ,
\quad
\qcc_7 := -\qcc_5,
\quad
\qcc_8 := -\qcc_6 .
%  &= \bigl((1,1),(-1,1),(-1,-1),(1,-1)\bigr).
\end{aligned}
\label{eq:d2q9-velocities}
\end{equation}
The D2Q9 discrete velocity set is illustrated in
\cref{fig:d2q9-velocities}. In the so-called \emph{lattice units}, $\dxLB=\dtLB=1$, so that
$\clattice=1$ for convenience.

%and the speed of sound is
%$\cs=\clattice/\sqrt{3}=1/\sqrt{3}$.

A complete LBE update consists of two key successive 
ingredients ---collision and streaming,
\begin{subequations}
\label{eq:lbe-evolution}
\begin{align}
	&
  \falpha^*(\vrj,\, t_n)
  = \falpha(\vrj,\, t_n)
     + \dtLB\Omega_\alpha(\qf(\vrj,\, t_n)),
\label{eq:lbe-collision}\\
  & \falpha(\vrj+\qxi_\alpha\dtLB,t_{n+1})
  = \falpha^*(\vrj,\, t_n),
  \quad
  \forall \vrj \in \dxLB\mathbb{Z}_2,
  \ \vrj +\qxi_\alpha\dtLB \in \dxLB\mathbb{Z}_2,
\label{eq:lbe-streaming}
\end{align}
\end{subequations}
where $\Omega_\alpha$ is the collision operator and the
superscript $*$ denotes the post-collision state.
Collision modifies the populations locally, while streaming
transports the post-collision distributions to neighboring grid
nodes along their respective discrete velocities.
Boundary conditions and other possible interactions
complete the update through prescribed procedures \cite{lallemand2021lattice}.
%by prescribing the
%missing incoming populations.

The moment description provides a general framework for constructing
the collision operator in the LBE.
This approach goes back to Maxwell's treatment of gas dynamics
in terms of moments of the molecular velocity distribution
\cite{Maxwell1867}.
The lowest-order moments determine the hydrodynamic variables,
while higher-order moments describe fluxes and other kinetic degrees
of freedom.
%characterize the remaining kinetic degrees of freedom.
Collision preserves the conserved moments and relaxes the
nonconserved ones.
%Different collision models employ different moment representations.
The multiple-relaxation-time (MRT) model uses an invertible linear
transformation $\tensorfont{M}$ to map distributions to their moments:
$\qvector{m}=\tensorfont{M}\qf$ and
$\qf=\tensorfont{M}^{-1}\qvector{m}$
\cite{Lallemand2000pre,dHumieres2002ptrsla}.
The relaxation frequencies of the nonconserved moments determine
transport coefficients and the damping of higher-order kinetic modes.
%where $\tensorfont{M}$ is the moment transformation matrix, 
%and assigns relaxation rates to the
%nonconserved moments
%The central-moment model takes moments about the local fluid
%velocity \cite{Geier2006pre}.
The two-relaxation-time (TRT) model separates the moments into even- and
odd-order ones
%populations into
%even and odd parts under velocity reversal
\cite{Ginzburg2008cicp427}.
The cumulant model uses cumulants of the velocity distribution
instead of moments \cite{Geier2015cma}.
These models differ in their collision operators, while retaining
the basic collision-streaming procedure of the LBE.

The present work considers 
%TRT collision with 
the preconditioned LBE (PLBE) \cite{guo2004preconditioned}
with the TRT collision model, which is a special case of the MRT model.
The macroscopic conserved variables in the system are mass density $\rho$ or
the pressure $\phat$ and the flow velocity $\qu$
\begin{equation}
  \phat = \cs^2 \sum_{\alpha=0}^{8}\falpha,
  \qquad
  \qu=\sum_{\alpha=0}^{8}\qxi_\alpha \falpha ,
\label{eq:trt-plbe-moments}
\end{equation}
%
%where $\phat$ is the pressure variable appearing in the
%equilibrium and $\qu$ is the macroscopic velocity.
where $\cs$ is the speed of sound in the system.
The equilibrium distribution in the PLBE is:
\begin{equation}
  \fazero(\phat,\qu;\gammaPLBE)
  =\frac{\walpha}{\cs^2}
   \left[
      \phat+\qxi_\alpha\cdot\qu
      +\frac{1}{2\gammaPLBE}
       \left(
          \frac{(\qxi_\alpha\cdot\qu)^2}{\cs^2}
          -\qu\cdot\qu
       \right)
   \right],
\label{eq:trt-plbe-equilibrium}
\end{equation}
where $\gammaPLBE \in (0,\, 1]$ is the kinetic preconditioning parameter.
%$0<\gammaPLBE\leq1$.
The equilibrium for the incompressible LBE \cite{HeXY1997jsp88} is recovered
with $\gammaPLBE=1$. The weight coefficients $\{w_\alpha\}$ are: 
$w_0=4/9$, $\walpha=1/9$ for $1\leq\alpha\leq4$, and
$\walpha=1/36$ for $5\leq\alpha\leq8$.

The discrete velocity set for the LBE, 
$\mathbb{X}_q := \{ \qxi_\alpha | 0 \leq \alpha \leq q-1 \}$
is symmetric, \emph{i.e.},
$\mathbb{X}_q = -\mathbb{X}_q$.
For convenience, denote $\bar{\alpha}$ the direction
opposite to $\alpha$, such that $\qcc_{\bar{\alpha}} := - \qcc_\alpha$. 
%To define TRT collision, let $\bar\alpha$ denote the direction
%opposite to $\alpha$, such that $\qcc_{\bar\alpha}=-\qcc_\alpha$.
The distribution functions may be split into even and odd parts,
corresponding to their contributions to even- and odd-order moments,
respectively,
%We split each population and its equilibrium into even and
%odd parts with respect to velocity reversal
\cite{Ginzburg2008cicp427}:
\begin{equation}
%\begin{aligned}
  f_\alpha^{\pm} = \frac{\falpha \pm \fbaralpha}{2},
  \quad
  \faeqpm =\frac{\fazero \pm \fbaralpha^{\ssstyle(0)}}{2},
%  \fminus&=\frac{\falpha-\fbaralpha}{2},
%  &\feqminus&=\frac{\fazero-\fbaralpha^{\ssstyle(0)}}{2}.
%\end{aligned}
\label{eq:trt-parity}
\end{equation}
%
%Here the superscript $(0)$ denotes equilibrium.
For the TRT collision model, the post-collision distribution is:
\begin{equation}
  \falpha^*
  = \falpha
  -\omegaplus \bigl( \fplus-\feqplus \bigr)
  -\omegaminus \bigl( \fminus-\feqminus \bigr),
  \label{eq:trt-collision}
\end{equation}
where $\omegaplus$ and $\omegaminus$ are the relaxation
frequencies for the even and odd order moments, respectively.
The distribution for the particles at rest, $f_0$,
for which $\bar\alpha=\alpha=0$,
has no odd part.
When $\omegaplus=\omegaminus$, \cref{eq:trt-collision}
reduces to the collision model with a single-relaxation-time (SRT) or the 
lattice BGK model.

%collision with the same equilibrium.

For the D2Q9 LBE model, the viscosity in
the lattice units is
\begin{equation}
  \viscLB
  =\frac{\gammaPLBE}{3}
    \left(\frac{1}{\omegaplus}-\frac{1}{2}\right) ,
  \label{eq:trt-lattice-viscosity}
\end{equation}
and the speed of sound is:
\begin{equation}
	\cs = \frac{1}{\sqrt{3}} c .
\end{equation}

% For a square cavity resolved with $\Ncavity\times\Ncavity$
% cells and lid speed $\Ulid$ in lattice units, the Reynolds
% number and the grid Reynolds number are
% \begin{equation}
%   \Rey=\frac{\Ulid \Ncavity}{\viscLB},
%   \qquad
%   \Reygrid:=\frac{\Rey}{\Ncavity}.
%   \label{eq:trt-grid-reynolds}
% \end{equation}

The collision, streaming, boundary conditions, and possible 
forcing effects together define a nonlinear discrete mapping $\Tdelta$
for the LBE:
\begin{equation}
  \qf^{m+1}=\Tdelta(\qf^{m}),
  \label{eq:trt-discrete-map}
\end{equation}
where $m$ denotes the successive discrete time $t_m$.

\subsection{Residual equation and nonlinear smoothing}
\label{subsec:trt-residual}

We define the unscaled steady-state residual directly from the complete
update mapping:
\begin{equation}
  \Rmap(\qf):=\qf-\Tdelta(\qf)
  := \dtLB\Adelta(\qf).
\label{eq:trt-unscaled-residual}
\end{equation}
Thus a steady discrete state satisfies $\Rmap(\qf)=\qvector{0}$. 
The weighted nonlinear iteration uses the same residual:
\begin{equation}
  \qf^{k+1}
  =\qf^{k}-\weightLBE\Rmap(\qf^{k}),
  \label{eq:trt-weighted-smoother}
\end{equation}
where $k$ is a nonlinear iteration index and $\weightLBE$ is the
smoother weight. For $\weightLBE \ne 1$, this is a
weighted fixed-point iteration, not physical-time stepping.

\subsection{FAS multigrid with fixed relaxation frequencies}
\label{subsec:trt-fas}

We use a hierarchy of uniform, cell-centered Cartesian grids with
factor-two coarsening.  The FAS construction follows the nonlinear
multigrid treatment of the LBE \cite{mavriplis2006multigrid,gsell2020multigrid,JiLuo2026cmame}.  Superscripts ``h" and ``H" denote two adjacent
fine and coarse FAS grids, respectively.  Each grid has its own
nonlinear TRT-PLBE update and residual operator.
The coarser operators are formed by rediscretization of the
collision--streaming map, including the boundary conditions, rather than
by a Galerkin construction of a Jacobian matrix.

The grid spacing and time step on successive FAS levels
follow acoustic scaling, \ie, $\dxH = 2 \dxh$ and $\dtH=2\dth$.
The TRT relaxation frequencies $\omegaplus$ and $\omegaminus$,
as well as the preconditioning parameter $\gammaPLBE$,
are kept at their finest-grid values on all FAS levels. The resulting 
coarse-grid problem is an auxiliary defect equation,
not an independent solution with the same physical Reynolds number.

The defect $\FASdefect^{\mathrm{h}}$ denotes the right-hand side on a fine FAS level. Thus, the corresponding nonlinear equation is
$\Rmap^{\mathrm{h}}(\qf^{\mathrm{h}})=\FASdefect^{\mathrm{h}}$. On the finest grid the defect is zero. For a current fine-grid iterate $\qf^{\mathrm{h}}$, the coarse-grid defect is constructed as
\begin{equation}
  \FASdefect^{\mathrm{H}}
   =\Rmap^{\mathrm{H}}(\Phistate\qf^{\mathrm{h}})
    -\frac{\dtH}{\dth}\Phires
       \bigl[\Rmap^{\mathrm{h}}(\qf^{\mathrm{h}})].
\label{eq:trt-fas-defect}
\end{equation}
Here $\Phistate$ restricts the distribution-function state, and
$\Phires$ restricts the nonlinear residual. The factor $\dtH/\dth=2$ follows from
$\Rmap=\dtLB\Adelta$ for the collision--streaming equations.
The weighted smoothing iteration on any FAS level is
\begin{equation}
  \qf^{\mathrm{h},k+1}
    =\qf^{\mathrm{h},k}
     -\weightLBE\bigl[\Rmap^{\mathrm{h}}(\qf^{\mathrm{h}, k})
                          -\FASdefect^{\mathrm{h}}\bigr].
  \label{eq:trt-fas-smoothing}
\end{equation}

After pre-smoothing, the coarse-grid initial state is
$\qf^{\mathrm{H},0}=\Phistate\qf^{\mathrm{h},\mathrm{pre}}$.
The coarse equation with the right-hand side of
\cref{eq:trt-fas-defect} is then solved approximately.
Its correction is returned to the fine grid by
\begin{equation}
  \widetilde{\qf}^{\mathrm{h}}
   =\qf^{\mathrm{h},\mathrm{pre}}
    +\Pcor\bigl(\qf^{\mathrm{H}}-\qf^{\mathrm{H}, 0}\bigr),
\label{eq:trt-coarse-correction}
\end{equation}
where $\Pcor$ is the correction prolongation operator.
Post-smoothing is applied to $\widetilde{\qf}^{\mathrm{h}}$ using
\cref{eq:trt-fas-smoothing}.  The coarsest level is treated by a
prescribed number of weighted nonlinear iterations rather than an
exact nonlinear solve.

\begin{figure}[tbp]
\centering
\resizebox{\linewidth}{!}{%
\begin{tikzpicture}[
  levelguide/.style={
    gray!55, densely dotted, line width=0.35pt
  },
  descending/.style={
    draw=black, line width=0.95pt,
    line cap=round
  },
  ascending/.style={
    draw=black, line width=0.95pt,
    dash pattern=on 3.2pt off 1.6pt,
    line cap=round
  },
  stagelabel/.style={font=\scriptsize},
  subtitle/.style={font=\small},
  coarsepoint/.style={
    fill=black, inner sep=1.4pt, circle
  }
]

% ---------- (a) V-cycle ----------
\begin{scope}[xshift=0cm]
  \node[subtitle] at (2.05,0.99) {(a) V-cycle};

  \foreach \y/\name in {0/L_1,-1/L_2,-2/L_3} {
    \draw[levelguide] (0.1,\y)--(4.1,\y);
    \node[stagelabel,anchor=east]
      at (0,\y) {$\name$};
  }

  \coordinate (v0) at (.45,0);
  \coordinate (v1) at (1.25,-1);
  \coordinate (v2) at (2.05,-2);
  \coordinate (v3) at (2.85,-1);
  \coordinate (v4) at (3.65,0);

  \draw[descending] (v0)--(v1)--(v2);
  \draw[ascending]  (v2)--(v3)--(v4);

  \node[coarsepoint] at (v2) {};

  \node[stagelabel,above=3pt]
    at (v0) {pre-smoothing};
  \node[stagelabel,above=3pt]
    at (v4) {post-smoothing};
\end{scope}

% ---------- (b) W-cycle ----------
\begin{scope}[xshift=5.2cm]
  \node[subtitle] at (2.45,.99) {(b) W-cycle};

  \foreach \y/\name in {0/L_1,-1/L_2,-2/L_3} {
    \draw[levelguide] (0.1,\y)--(4.95,\y);
    \node[stagelabel,anchor=east]
      at (0,\y) {$\name$};
  }

  \coordinate (w0)  at (.22,0);
  \coordinate (w1)  at (.60,-1);
  \coordinate (w2)  at (.98,-2);
  \coordinate (w3)  at (1.36,-1);
  \coordinate (w4)  at (1.74,-2);
  \coordinate (w5)  at (2.12,-1);
  \coordinate (w6)  at (2.50,0);
  \coordinate (w7)  at (2.88,-1);
  \coordinate (w8)  at (3.26,-2);
  \coordinate (w9)  at (3.64,-1);
  \coordinate (w10) at (4.02,-2);
  \coordinate (w11) at (4.40,-1);
  \coordinate (w12) at (4.78,0);

% W-cycle: restriction
\draw[descending] (w0)--(w1)--(w2);
\draw[descending] (w3)--(w4);
\draw[descending] (w6)--(w7)--(w8);
\draw[descending] (w9)--(w10);

% W-cycle: prolongation
\draw[ascending] (w2)--(w3);
\draw[ascending] (w4)--(w5)--(w6);
\draw[ascending] (w8)--(w9);
\draw[ascending] (w10)--(w11)--(w12);

% Coarsest-grid solves
\node[coarsepoint] at (w2) {};
\node[coarsepoint] at (w4) {};
\node[coarsepoint] at (w8) {};
\node[coarsepoint] at (w10) {};
\end{scope}

% ---------- (c) F-cycle ----------
\begin{scope}[xshift=11.25cm]
  \node[subtitle] at (2.45,.99) {(c) F-cycle};

  \foreach \y/\name in {0/L_1,-1/L_2,-2/L_3} {
    \draw[levelguide] (0.1,\y)--(4.95,\y);
    \node[stagelabel,anchor=east]
      at (0,\y) {$\name$};
  }

  \coordinate (f0)  at (.22,0);
  \coordinate (f1)  at (.675,-1);
  \coordinate (f2)  at (1.13,-2);
  \coordinate (f3)  at (1.585,-1);
  \coordinate (f4)  at (2.04,-2);
  \coordinate (f5)  at (2.495,-1);
  \coordinate (f6)  at (2.95,0);
  \coordinate (f7)  at (3.405,-1);
  \coordinate (f8)  at (3.86,-2);
  \coordinate (f9)  at (4.315,-1);
  \coordinate (f10) at (4.77,0);

% F-cycle: restriction
\draw[descending] (f0)--(f1)--(f2);
\draw[descending] (f3)--(f4);
\draw[descending] (f6)--(f7)--(f8);

% F-cycle: prolongation
\draw[ascending] (f2)--(f3);
\draw[ascending] (f4)--(f5)--(f6);
\draw[ascending] (f8)--(f9)--(f10);

% Coarsest-grid solves
\node[coarsepoint] at (f2) {};
\node[coarsepoint] at (f4) {};
\node[coarsepoint] at (f8) {};
\end{scope}

% ---------- Common legend ----------
\draw[descending] (1.3,-3.0)--(2.0,-3.0);
\node[stagelabel,anchor=west]
  at (2.12,-3.0) {restriction (downward)};

\draw[ascending] (6.0,-3.0)--(6.7,-3.0);
\node[stagelabel,anchor=west]
  at (6.82,-3.0) {prolongation (upward)};

\node[coarsepoint] at (11.0,-3.0) {};
\node[stagelabel,anchor=west]
  at (11.3,-3.0) {coarse solve};

\end{tikzpicture}%
}

\caption{Schematic traversal patterns of three-level
FAS multigrid cycles: (a) V-cycle, (b) W-cycle, and
(c) F-cycle.}
\label{fig:mg_cycles}
\end{figure}
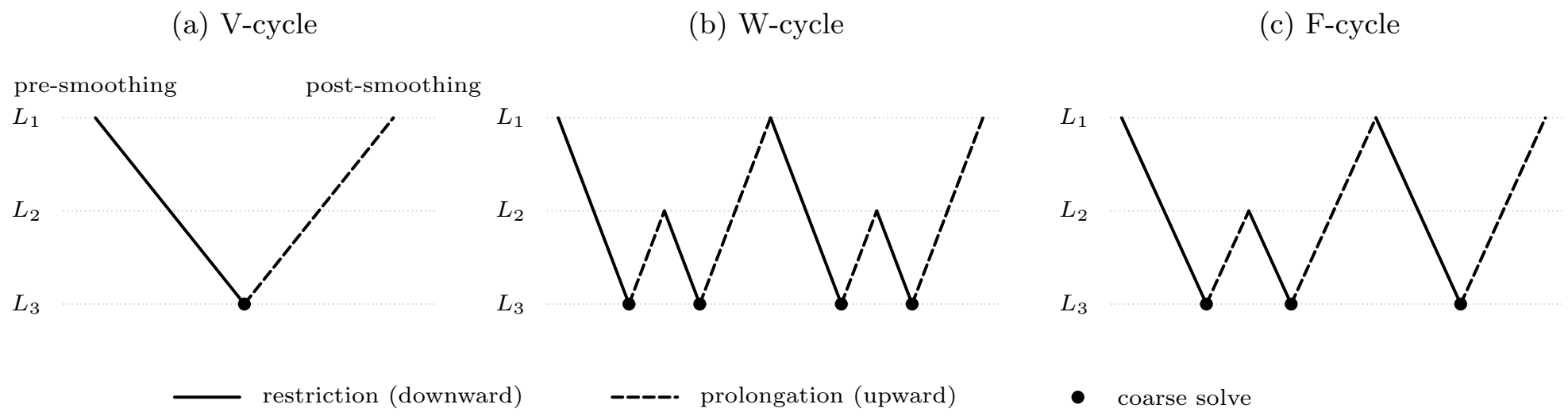

The FAS procedure described above corresponds to a V-cycle.
Other cycle types, including W-cycles and F-cycles, can also
be used. These cycles differ in their coarse-grid traversal
patterns, as illustrated in \cref{fig:mg_cycles}.
All three cycle types are considered in the numerical experiments.

\subsection{Cell-centered transfer and boundary treatment}
\label{subsec:trt-transfer}

The state restriction, residual restriction, and correction
prolongation act directly on each distribution function $\falpha$.
There is no reconstruction or rescaling of the nonequilibrium part
during these transfers.

For a fine-to-coarse transfer in two dimensions, the state restriction
uses a four-cell average.  The cell-centered construction follows
Mohr and Wienands \cite{mohr2004cell}:
\begin{equation}
  \Phistate=\frac{1}{4}
  \begin{bmatrix}1&1\\1&1\end{bmatrix}.
  \label{eq:trt-state-restriction}
\end{equation}
At interior coarse cells, the residual restriction uses the
bilinear stencil
\begin{equation}
  \Phires=\frac{1}{64}
  \begin{bmatrix}
     1&3&3&1\\
     3&9&9&3\\
     3&9&9&3\\
     1&3&3&1
  \end{bmatrix}.
  \label{eq:trt-residual-restriction}
\end{equation}
The matrix passed to the FAS residual restriction incorporates the
factor $\dtH/\dth=2$ in \cref{eq:trt-fas-defect}.
The correction prolongation $\Pcor$ uses the cell-centered
bilinear interpolation described in \cite{mohr2004cell}.

At a physical boundary, the restriction stencil of $\Phires$ is modified to
preserve constants. For example, the first one-dimensional
coarse-grid restriction row has weights $(1/2,3/8,1/8)$ on the
first three fine cells. The analogous rule is applied at the
opposite boundary and by tensor products in two dimensions.
Boundary coarse-grid corrections outside
the computational domain are set to zero.

\subsection{Relaxation parameters and computational procedure}
\label{subsec:trt-parameters}
The TRT relaxation frequencies, $\omega^+$ and $\omega^-$, affect 
both the physical viscosity and the nonlinear smoothing behavior. 
For a characteristic velocity $\Uzero$ in lattice units and a characteristic 
length resolved by $\Ncavity$ cells on the finest grid, we define the Reynolds number
$\Rey$ and the grid Reynolds number $\Rey^*$ as
\begin{equation}
  \Rey=\frac{\Uzero\Ncavity}{\viscLB},
  \qquad
  \Reygrid:=\frac{\Rey}{\Ncavity}
           =\frac{\Uzero}{\viscLB}.
\label{eq:trt-grid-reynolds}
\end{equation}
At a prescribed grid Reynolds number, $\Uzero$, and the
preconditioning parameter $\gammaPLBE$ provide two ways to adjust
$\omegaplus$ to improve the smoothing properties of
\cref{eq:trt-fas-smoothing}.
% The TRT relaxation frequencies affect both the physical viscosity
% and the nonlinear smoothing behavior. The characteristic speed $\Uzero$ and preconditioning parameter $\gammaPLBE$ provide two ways to adjust $\omegaplus$ at a prescribed grid Reynolds number to improve the smoothing property of \cref{eq:trt-fas-smoothing}. 
And $\omegaminus$ can be chosen independently to eliminate the local amplification of odd nonequilibrium components. \cref{eq:trt-lattice-viscosity} and
\cref{eq:trt-grid-reynolds} give
\begin{equation}
  \omegaplus
   =\left(\frac{1}{2}
        +\frac{3\Uzero}{\gammaPLBE\Reygrid}\right)^{-1}.
\label{eq:trt-omega-selection}
\end{equation}
For a fixed $\gamma$,
decreasing $\Uzero$
%with $\gammaPLBE$ fixed 
increases $\omegaplus$; whereas 
for a fixed $\Uzero$,
decreasing $\gammaPLBE$ decreases $\omegaplus$.

For fixed $\Uzero$ and $\gammaPLBE$, $\omegaplus$ approaches $2$ as the
grid Reynolds number $\Reygrid$ increases, whereas it decreases
as $\Reygrid$ becomes smaller, according to \cref{eq:trt-omega-selection}.
In both limits, the smoothing
property of the LBE iteration may deteriorate. A small
$\omegaplus$ leads to weak relaxation of the even-order moments,
%nonequilibrium components, 
while $\omegaplus$ close to $2$ causes severe
%strong
over-relaxation. To improve the smoothing property, we use
two different adjustments. For large $\Reygrid$, we reduce
the preconditioning parameter $\gammaPLBE$ to decrease
$\omegaplus$. For small $\Reygrid$, we instead reduce the
characteristic velocity $\Uzero$, while
keeping $\gammaPLBE=1$, to increase $\omegaplus$. In both
cases, the viscosity is adjusted to maintain the prescribed
$\Rey$. These adjustments allow $\omegaplus$ to remain in
a range favorable for nonlinear smoothing. The present
numerical tests suggest that $\omegaplus \in (1.53,\, 1.60)$
%$1.53$--$1.60$ gives 
leads to favorable FAS convergence \cite{JiLuo2026cmame}.

The odd relaxation frequency $\omegaminus$ can be selected
independently of $\omegaplus$. For the weighted LBE smoother
with the weight $\weightLBE$, the local amplification factor
of the nonconserved odd-order moments, with fixed equilibrium
moments, is $(1-\weightLBE\omegaminus)$. Therefore this
amplification 
factor may be nullified with $\weightLBE\omegaminus=1$. 
We therefore choose $\omegaminus=1/\weightLBE$, which gives
$\omegaminus=1.25$ for $\weightLBE=0.8$. This choice provides
complete local relaxation of the odd-order nonconserved moments
%nonequilibrium component
under the aforementioned assumption.

For a given grid Reynolds number $\Reygrid$ and a chosen
target $\omegaplus$ in the admissible interval $(1.53,\, 1.60)$,
these two adjustments can be written, respectively, as
\begin{equation}
    \Uzero=\frac{\gammaPLBE\Reygrid}{3}\left(\frac{1}{\omegaplus}-\frac{1}{2}\right), \qquad
    \gammaPLBE=\frac{3\Uzero}{\Reygrid\left(1/\omegaplus-1/2\right)}.
  \label{eq:trt-target-selection}
\end{equation}
The first relation applies when $\gammaPLBE$ is prescribed, and the
second when $\Uzero$ is prescribed.  The target is feasible only if
\cref{eq:trt-gamma-bound} is satisfied.

In the configurations considered in what follows, the odd relaxation frequency
is tied to the nonlinear smoother weight by
\begin{equation}
  \omegaminus=\frac{1}{\weightLBE},
  \qquad \weightLBE=0.8,
  \label{eq:trt-odd-weight-choice}
\end{equation}
giving $\omegaminus=1.25$.
The relation specifies the implemented smoother, not an analytical
optimality condition. 

For kinetic preconditioning, we monitor the estimated preconditioned
Mach number
\begin{equation}
  \Mapre:=\frac{\Uzero}{\cs\sqrt{\gammaPLBE}}
  < 0.3 .
%  \qquad \Mapre<0.3.
\label{eq:trt-mach-bound}
\end{equation}
This condition restricts the admissible values of
$\gammaPLBE$ as follows:
\begin{equation}
  \gammaPLBE>
  \left(\frac{\Uzero}{0.3\cs}\right)^2.
\label{eq:trt-gamma-bound}
\end{equation}
The bound is a parameter-selection condition rather than a 
result based on analysis.
%of accuracy or nonlinear stability.

The cycles are implemented in PETSc SNESFAS
\cite{brune2015composing}.  In the FAS cycle,
$\smoothpre$ and $\smoothpost$ count pre- and post-smoothing steps,
and $\smoothcoarse$ counts coarsest-grid steps. Every level uses the 
same $\weightLBE$.  The finest-grid solve is stopped by the 
Euclidean norm of its distribution function residual \cref{eq:trt-unscaled-residual},
with relative tolerance $10^{-8}$ and absolute tolerance $10^{-12}$.

\section{Numerical results}
\label{sec:trt-numerical-results}

% We consider steady two-dimensional lid-driven cavity flow in the unit
% square. The lid moves with speed $\Ulid$ and the remaining walls are
% stationary. The reported velocities are normalized by the lid speed,
% and spatial coordinates by the cavity side length. The computations
% use uniform cell-centered grids. 

We consider 
the two-dimensional lid-driven cavity flow as a benchmark
for solving the steady incompressible Navier--Stokes equations,
without assuming temporal stability at high Reynolds numbers.
The problem is a classical multigrid benchmark \cite{ghia1982high}.
While the boundary geometry is simple, the flow exhibits various 
complicated features, such as thin boundary layers and multiple
flow separation vortices. 
%combining simple geometry with thin boundary layers and corner
%vortices at high Reynolds numbers.
The flow domain is a square with no-slip walls on all sides.
The top side moves with speed
$\Ulid$ and the remaining walls are stationary. 
%The reported velocities are normalized by the lid speed, 
%and spatial coordinates by the cavity side length. 
The computations use uniform cell-centered grids.

\subsection{Nonlinear convergence and computational efficiency}
\label{subsec:trt-high-re-convergence}

\begin{figure}[htbp]
\centering
\includegraphics[width=0.75\linewidth]{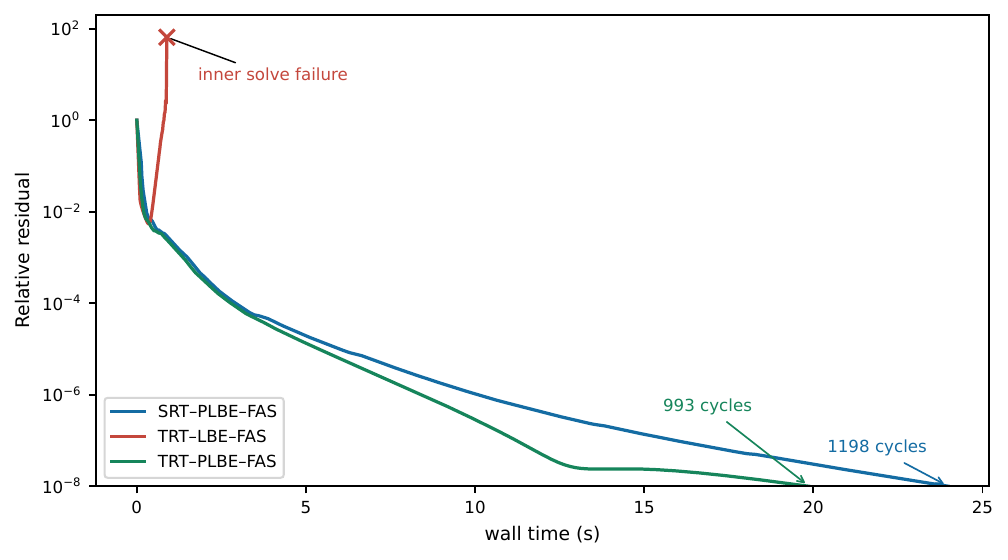}
\caption{2D lid-driven cavity flow at $\Rey=5000$, grid size $512^2$.
Convergence histories of three FAS smoother configurations 
with $\Ulid=0.05$. The TRT-LBE-FAS curve ends at the last
monitored residual before the inner solve failure.}
\label{fig:trt-eff-e1}
\end{figure}

We first examine the effect of the nonlinear smoother at
$\Rey=5000$ on a $512^2$ grid. All three FAS runs in
\cref{fig:trt-eff-e1} use six grid levels, a V-cycle with
$\smoothpre=\smoothpost=3$, and $\weightLBE=0.8$. 
For the TRT-PLBE smoother in FAS, we set $\omegaplus=1.6$, 
$\omegaminus=1.25$, and $\gamma=0.12288$.
The SRT-PLBE-FAS and TRT-PLBE-FAS configurations reach the prescribed
relative residual tolerance in 1198 and 993 cycles, respectively,
corresponding to wall-clock times of 23.99~s and 19.87~s. Thus, the 
TRT-PLBE-FAS reduces the solution time by about 17.16\% in this 
case.
The tested TRT-LBE-FAS configuration (without preconditioning) 
terminates with an inner
iteration failure and is not assigned a time to convergence.

\begin{figure}[htbp]
\centering
\includegraphics[width=0.75\linewidth]{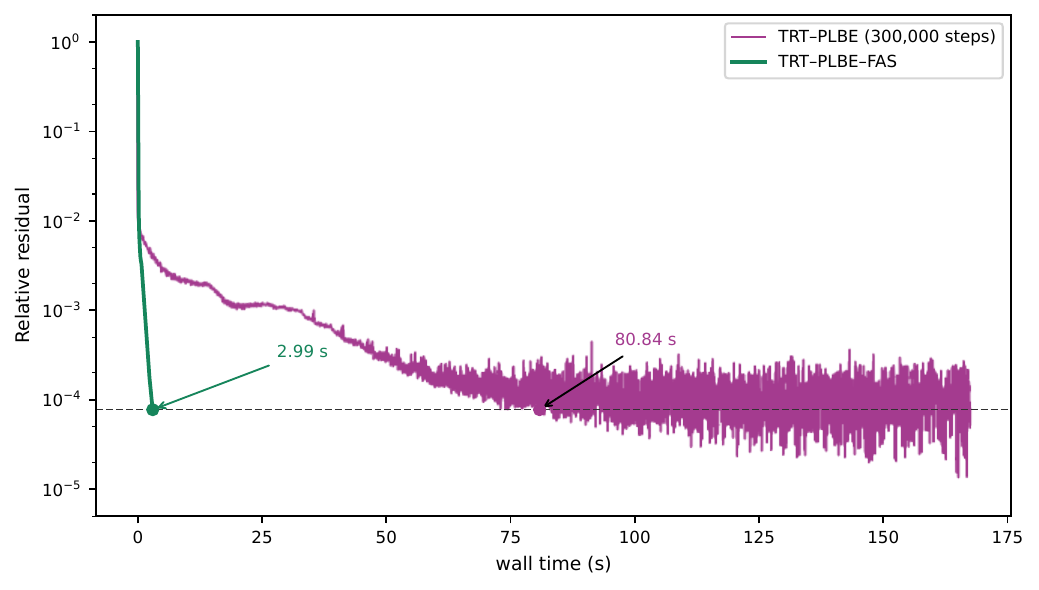}
\caption{Single-grid TRT-PLBE and TRT-PLBE-FAS at $\Rey=5000$
on a $512^2$ grid, with $\Ulid=0.05$. The dashed line marks
the common relative residual threshold
$7.77961\times10^{-5}$, defined by the final monitored residual
of the 300000-step single-grid run.}
\label{fig:trt-eff-e2}
\end{figure}

We next compare single-grid TRT-PLBE and TRT-PLBE-FAS at the same
$\Rey$ and grid resolution. The single-grid TRT-PLBE uses the same parameters
as the smoother. The TRT-PLBE does not reach the
$10^{-8}$ relative residual tolerance within 300000 steps
(167.37s), ending at $7.77961\times10^{-5}$.
Taking this final monitored residual as a common threshold, the
single-grid method first attains it after 145972 steps (80.839s),
while TRT-PLBE-FAS first attains it after 161 cycles (2.993s).
%The resulting first-passage time ratio is $27.01$.
So the speedup due to FAS is about 27 times.
In \cref{fig:trt-eff-e2}, the FAS curve is truncated at the first
crossing of this threshold, whereas the complete single-grid history
is shown. This ratio measures the times required to attain a common
residual level; it is not a speedup to the $10^{-8}$ stopping tolerance.

Unless otherwise stated,
%the below high-Reynolds-number cycle tests use fixed relaxation frequencies
the following parameters are used in subsequent tests:
$\omegaplus=1.6$, $\omegaminus=1.25$, preconditioning parameter
$\gammaPLBE=0.12288$, and smoother weight $\weightLBE=0.8$.
The runs summarized below use $\smoothpre=\smoothpost=5$ and $\smoothcoarse=40$.

The choice of FAS cycle becomes important in the high-Reynolds-number
regime. The V-cycle exhibits residual stagnation in some
of the high-Reynolds-number tests. The W-cycle
performs additional recursive coarse-grid
corrections, whereas the hybrid V/W-cycle applies
the more expensive recursion only on selected
levels. The F-cycle configuration employs
a distinct coarse-to-fine initialization followed
by V-cycle iterations.

\Cref{tab:trt-high-re-cycles} summarizes selected completed
runs. At $\Rey=10000$ on a $1024^2$ grid,
the F-cycle configuration reached the prescribed relative residual
tolerance of $10^{-8}$. At $\Rey=20000$ on a $2048^2$ grid, the
F-cycle configuration encountered a residual plateau, whereas both
the W-cycle and a hybrid V/W-cycle reached the same stopping
criterion. In the hybrid configuration, the intermediate FAS
levels 2--6 use W-cycles, level 1 uses a V-cycle, and level 0 retains
the coarsest-grid solver. In these recorded runs, the hybrid V/W required more cycles than W but less wall-clock time.

\begin{table}[htbp]
\centering
\footnotesize
\caption{Nonlinear convergence and computational efficiency for the
selected cavity-flow tests on one NVIDIA RTX 4090 GPU.}
\label{tab:trt-high-re-cycles}
\label{tab:trt-efficiency-summary}
\setlength{\tabcolsep}{4pt}
\begin{tabular}{@{}lrrrr@{}}
\toprule
Solver / FAS cycle & $\Rey$ & Grid
& Cycles or steps & Time (s) \\
\midrule
SRT-PLBE-FAS (V) & 5000 & $512^2$ & 1198 & 23.99 \\
TRT-LBE-FAS (V) & 5000 & $512^2$ & 45 & $0.90^{\dagger}$ \\
TRT-PLBE-FAS (V) & 5000 & $512^2$ & 993 & 19.87 \\
TRT-PLBE (single grid) & 5000 & $512^2$
& 300000 & $167.37^{\ddagger}$ \\
\midrule
TRT-PLBE-FAS (F-cycle) & 10000 & $1024^2$
& 481 & 83.29 \\
TRT-PLBE-FAS (F-cycle) & 20000 & $2048^2$
& -- & -- \\
TRT-PLBE-FAS (W-cycle) & 20000 & $2048^2$
& 121 & 180.52 \\
TRT-PLBE-FAS (hybrid V/W) & 20000 & $2048^2$
& 168 & 154.65 \\
\bottomrule
\end{tabular}
\par\smallskip
\begin{minipage}{\linewidth}\footnotesize
$^{\dagger}$Inner solve failure; the time is not a time to tolerance.
$^{\ddagger}$Step budget exhausted; the run did not reach the
target tolerance. The $\Rey=20000$ F-cycle run encountered a
residual plateau, and no convergence time is reported. The
single-grid count is in steps; all other counts are FAS cycles.
\end{minipage}
\end{table}

\begin{figure}[htbp]
\centering
\includegraphics[width=0.75\linewidth]{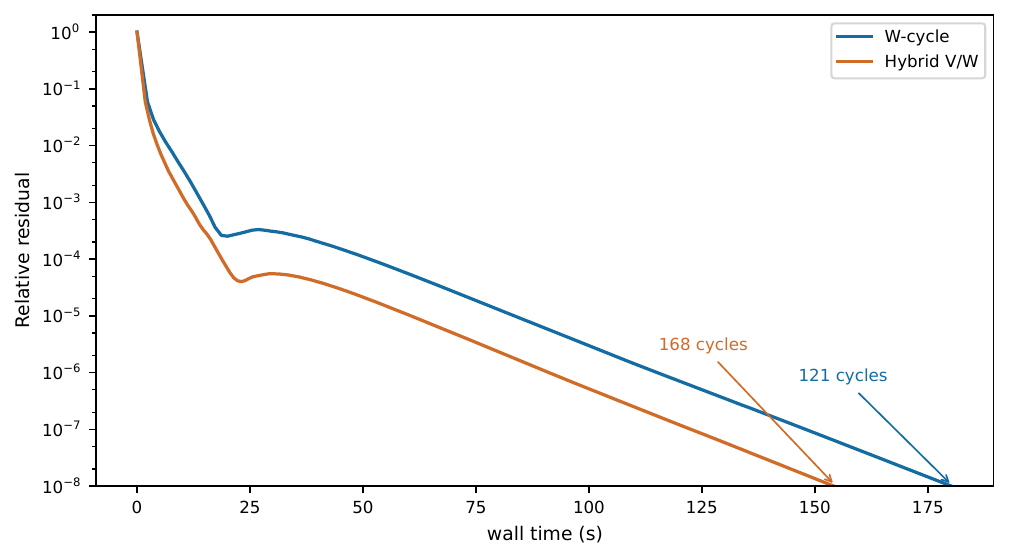}
\caption{Convergence histories for the
W-cycle and hybrid V/W configuration at $\Rey=20000$ on a
$2048^2$ grid, with $\Ulid=0.05$. The hybrid cycle reduces the measured total solve time by $14.33\%$.}
\label{fig:trt-eff-e3}
\end{figure}

For the paired $\Rey=20000$ tests, the W-cycle and hybrid V/W
configuration use the same executable and stopping tolerances.
As shown in \cref{fig:trt-eff-e3}, the W-cycle reaches the prescribed
tolerance in 121 cycles and 180.52~s. Hybrid V/W requires 168 cycles
but only 154.65~s, reducing total solve time by $14.33\%$ in this
single paired comparison. This is an end-to-end solve-time result;
it does not quantify the contributions of smoothing, transfer, and
coarse-grid correction, for which mutually exclusive timings are
not yet available.

To further examine the computational cost of different FAS
cycles, we measure one top-level FAS application. The computational time is decomposed into
smoothing, restriction, prolongation, coarse-grid solving,
residual evaluation, and other operations. To obtain
mutually exclusive stage timings, GPU synchronization is
introduced at the measurement boundaries.
\Cref{tab:trt-cycle-cost} reports the synchronized stage
times and, separately, the wall-clock times measured
without intermediate synchronization.

\begin{table}[htbp]
\centering
\scriptsize
\caption{Computational cost (ms) of one FAS
application at $\Rey=20000$ on a $2048^2$ grid.}
\label{tab:trt-cycle-cost}
\resizebox{\linewidth}{!}{%
\begin{tabular}{@{}lrrrrrrrr@{}}
\toprule
Cycle & Smooth & Restrict & Prolong & Coarse
& Residual & Other & Sync total & Wall time \\
\midrule
V-cycle
& 138.7 & 610.9 & 341.0 & 8.4
& 5.4 & 1.2 & 1105.6
& 1105 [1095, 1112] \\
W-cycle
& 472.1 & 622.3 & 351.4 & 1080.0
& 10.7 & 10.4 & 2546.9
& 2540 [2529, 2540] \\
Hybrid V/W
& 402.2 & 618.9 & 350.6 & 538.1
& 9.2 & 8.0 & 1927.1
& 1922 [1921, 1927] \\
\bottomrule
\end{tabular}%
}
\end{table}

The W-cycle is more expensive per application than
the V-cycle because of its additional recursive visits
to the coarse levels. The hybrid V/W-cycle reduces
this cost by limiting W-type recursion to selected
levels. In the synchronized measurements, the
difference between W and hybrid V/W is approximately
620ms per application, of which 542ms arises from
the coarsest-grid solver. Thus, the reduced number
of coarse-grid visits accounts for most of the
measured cost difference.

\subsection{Comparison with benchmark data}
\label{subsec:trt-cavity-accuracy}

\begin{figure}[htbp]
\centering
\begin{minipage}{0.98\linewidth}
\centering
\includegraphics[width=\linewidth]{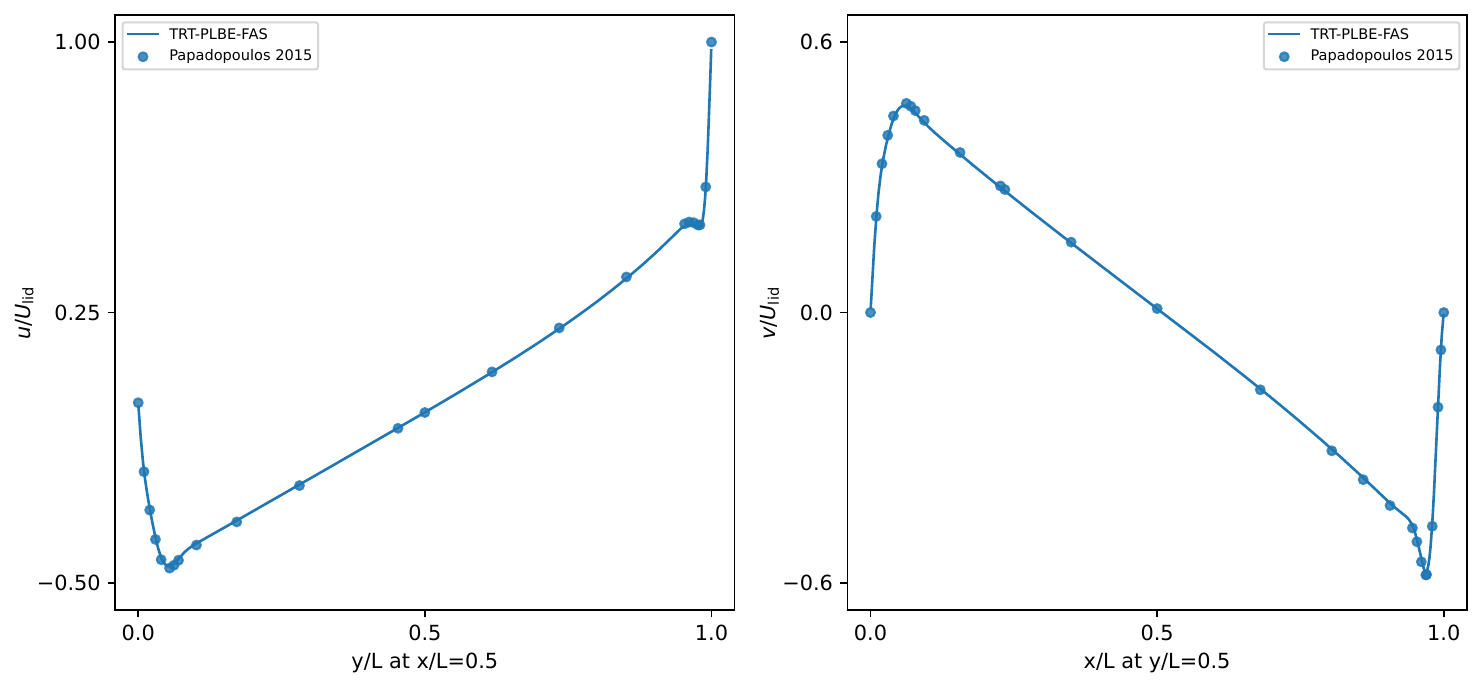}
(a) $\Rey=10000$, $1024^2$ grid
\end{minipage}

\begin{minipage}{0.98\linewidth}
    \centering
    \includegraphics[width=\linewidth]{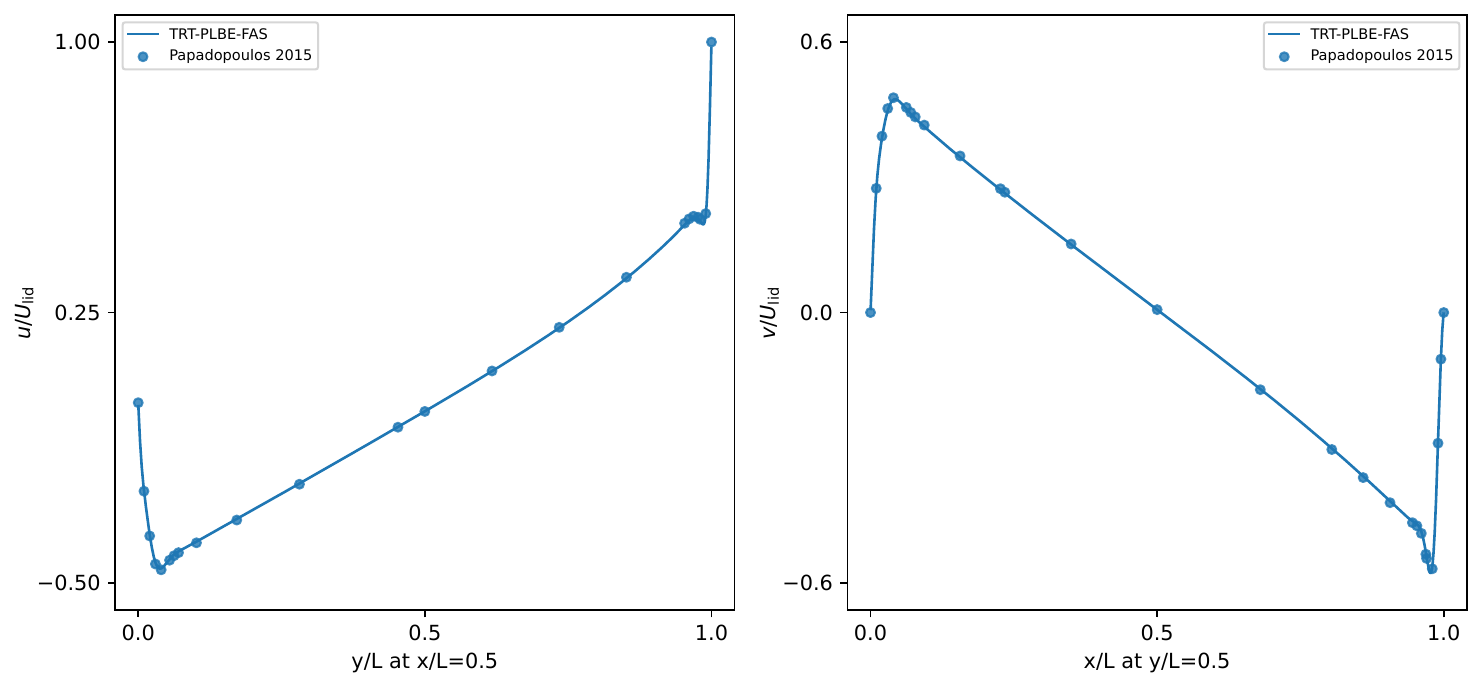}
(b) $\Rey=20000$, $2048^2$ grid
\end{minipage}
\setlength{\abovecaptionskip}{3pt}
\caption{2D lid-driven cavity flow at $\Rey=10000$ with grid size
$1024^2$ (top row)
and $\Rey=20000$ with grid size $2048^2$ (bottom row).
Normalized cavity centerline velocities.
(Left) $u(x/L=0.5,\, y/L)/\Ulid$ and (right)
$v(x/L,\, y/L=0.5)/\Ulid$.
%Each pair of panels shows $u/\Ulid$ along $x/L=1/2$ (left)
%and $v/\Ulid$ along $y/L=1/2$ (right), compared with
%the available reference data \cite{Papadopoulos2015}.
}
\label{fig:trt-centerline-comparison}
\end{figure}

We compare the converged solutions at $\Rey=10000$ (grid size $1024^2$,
F-cycle) and $\Rey=20000$ (grid size $2048^2$, hybrid cycle) with the
high-resolution cavity data of Papadopoulos
\cite{Papadopoulos2015}.
The reference data are consistent with the
high-order finite-element results of
Shapeev and Lin \cite{shapeev2009asymptotic},
which report closely matching vortex intensities
and locations at $\Rey=10000$ and $\Rey=20000$.
The centerline velocities are sampled at
$x/L=1/2$ and $y/L=1/2$. \Cref{fig:trt-centerline-comparison} shows the normalized profiles.
The maximum absolute discrepancies at the tabulated reference
sampling locations are $6.834\times10^{-3}$ for $u/\Ulid$ and
$7.497\times10^{-3}$ for $v/\Ulid$ at $\Rey=10000$, and
$5.029\times10^{-3}$ and $5.405\times10^{-3}$, respectively,
at $\Rey=20000$.

%

% The global streamfunction contours in
% \cref{fig:trt-global-streamlines} provide a comparison
% of flow structure. 
% % We define the streamfunction by
% % $u=\partial\psi/\partial y$ and
% % $v=-\partial\psi/\partial x$, with the common wall value set to
% % zero, and display $\psi/(\Ulid L)$. The streamfunction is
% % reconstructed from the velocity field by averaging two path
% % integrals. 
% At both Reynolds numbers, the principal recirculation
% occupies most of the cavity, accompanied by secondary corner
% eddies. And the quantitative comparisons are made with centerline velocities
% and streamfunction extrema.

\begin{figure}[htbp]
\centering
\includegraphics[height=0.375\linewidth]{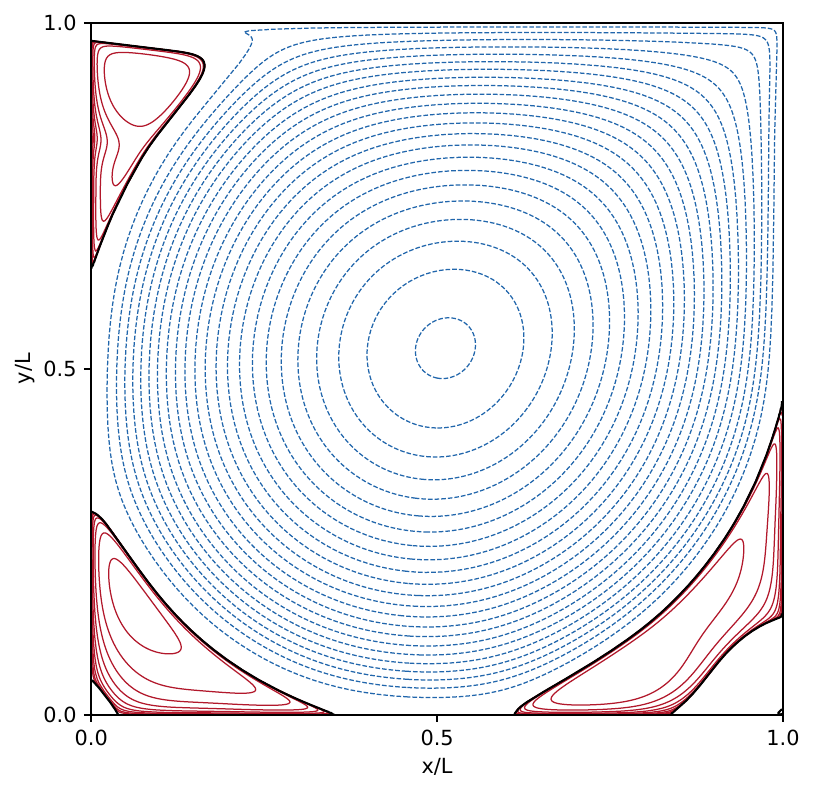}%
\includegraphics[height=0.375\linewidth]{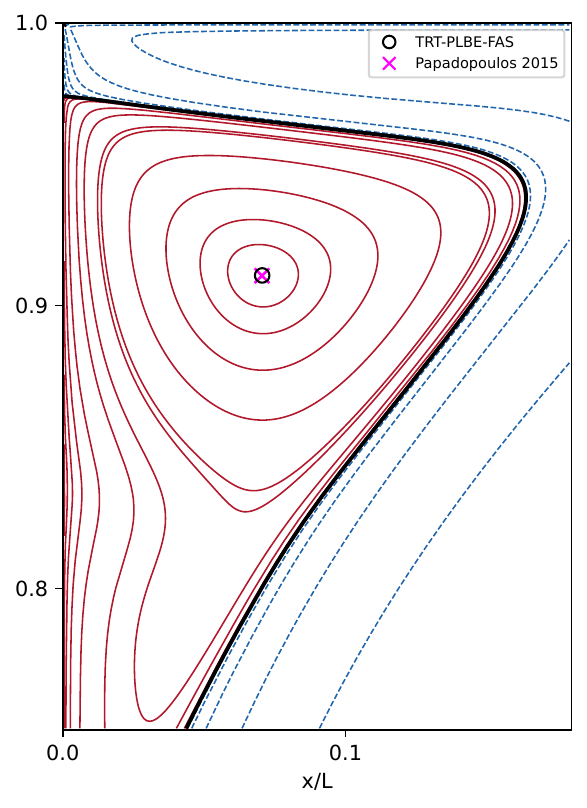}%
\includegraphics[height=0.375\linewidth]{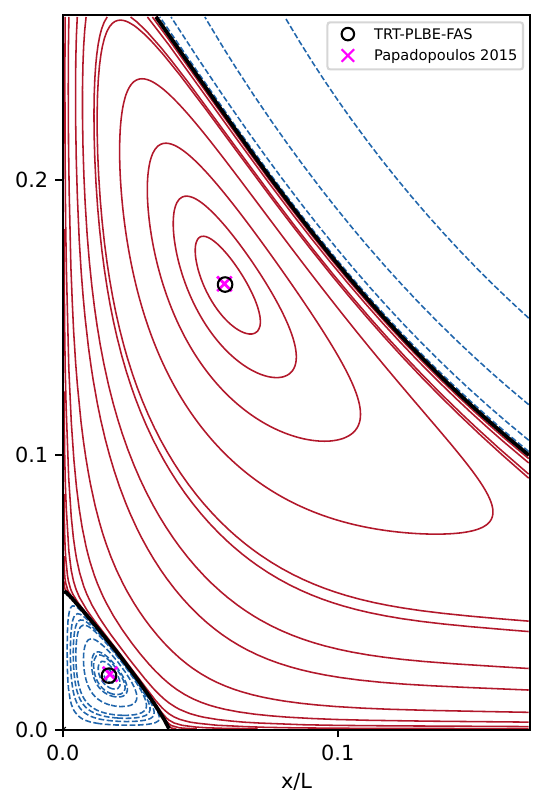}\\
\includegraphics[width=0.55\linewidth]{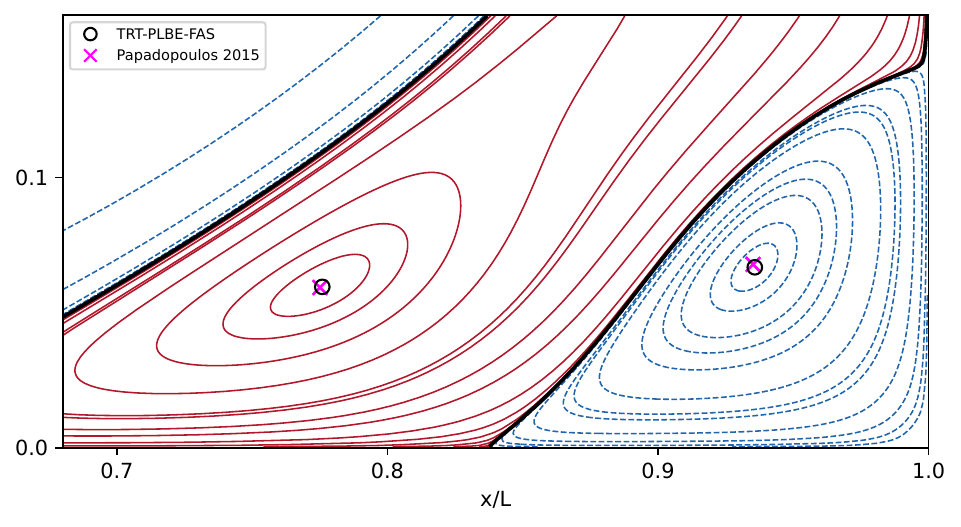}%
\caption{2D lid-driven cavity flow at $\Rey=10000$ with grid size
$1024^2$. The normalized streamfunction $\psi/(\Ulid L)$
contours in (from left to right, top to bottom):
(a) the entire flow domain, (b) upper-left corner,
(c) lower-left corner, and (d) lower-right corner.
Solid red and dashed blue contours
denote positive and negative values, respectively,
while thick black contours denote zero contours.
Open circles indicate computed vortex center positions, and crosses
indicate locations given in \cite{Papadopoulos2015}.
}
\label{fig:trt-streamlines-re10000}
\end{figure}

\begin{figure}[htbp]
\centering
\includegraphics[height=0.375\linewidth]{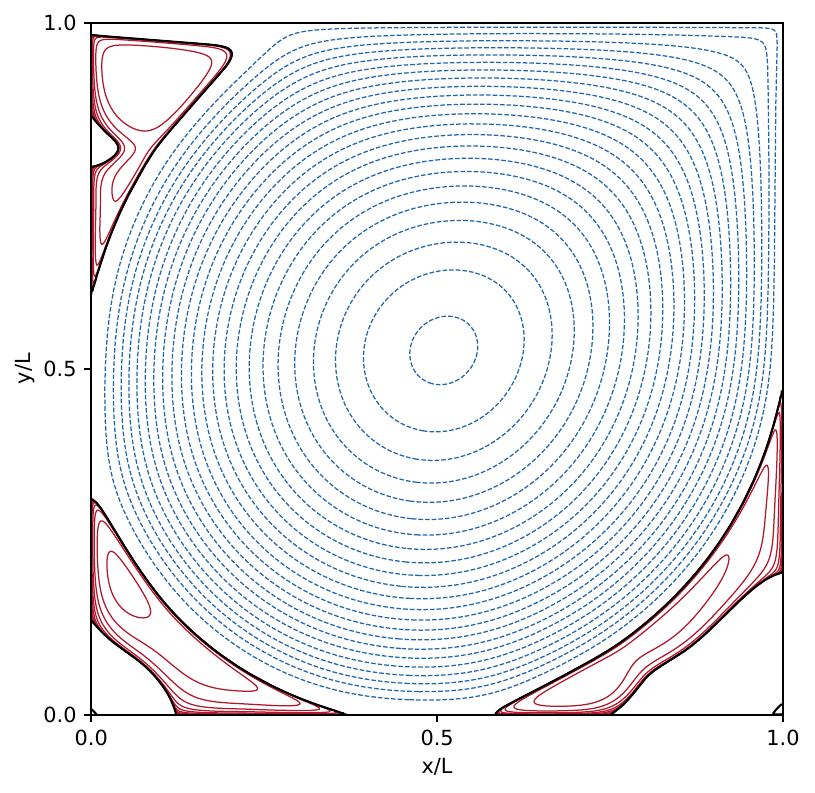}%
\includegraphics[height=0.375\linewidth]{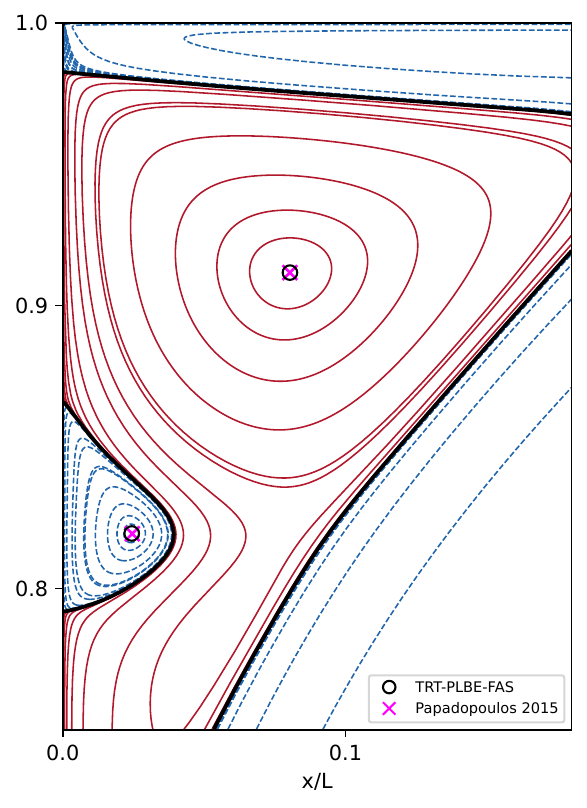}%
\includegraphics[height=0.375\linewidth]{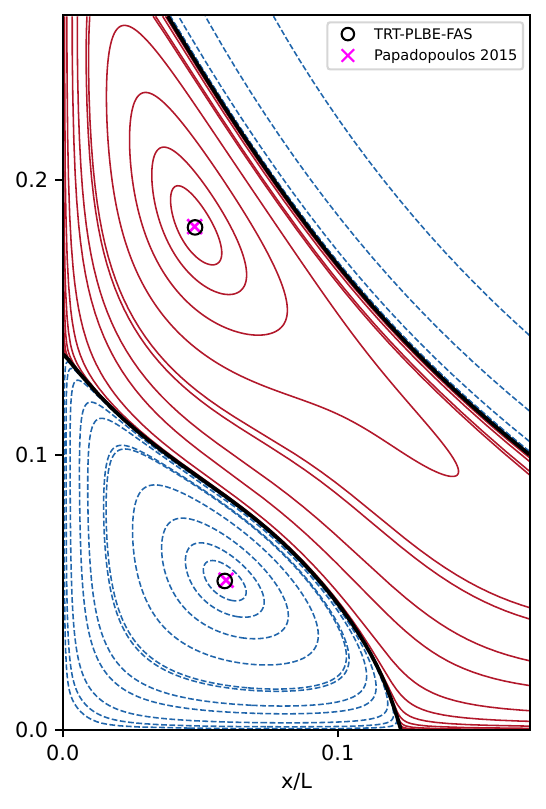}\\
\includegraphics[width=0.55\linewidth]{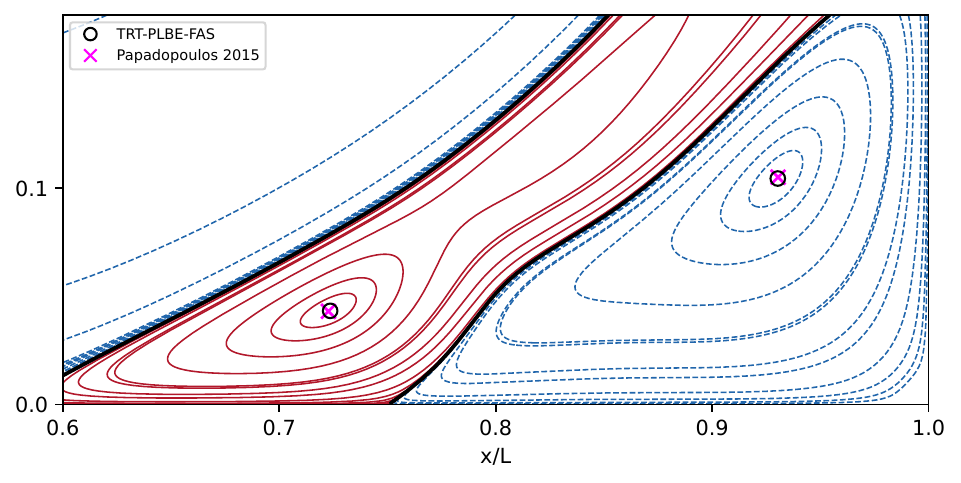}%
\caption{2D lid-driven cavity flow at $\Rey=20000$ with grid size
$2048^2$. Similar to \cref{fig:trt-streamlines-re10000}.}
%Streamfunction contours for lid-driven cavity flow
%at $\Rey=20000$ on a $2048^2$ grid:
%(a) global flow, (b) upper-left corner,
%(c) lower-left corner, and (d) lower-right corner.
%Contour styles and vortex-center markers follow those
%in \cref{fig:trt-streamlines-re10000}.}
\label{fig:trt-streamlines-re20000}
\end{figure}

The global and corner-region streamfunction contours at
$\Rey=10000$ and $20000$ are shown in
\cref{fig:trt-streamlines-re10000,fig:trt-streamlines-re20000},
respectively. At both Reynolds numbers, the primary
recirculation occupies most of the cavity, accompanied by
secondary corner vortices. The enlarged corner views show
the resolved vortex structures and compare the computed
vortex centers with the reference locations of Papadopoulos
\cite{Papadopoulos2015}.
Some of the smaller and weaker corner vortices reported
in the reference solution are not reliably resolved on
the present grids. Therefore, only reliably identified
vortex centers are marked and included in the
quantitative comparisons.

\Cref{tab:trt-vortex-extrema} quantifies the positions and
streamfunction extrema of the principal vortex and selected resolved
corner vortices at $\Rey=20000$. For the primary vortex,
the computed center $(0.50960,0.52669)$ differs from the reference
location $(0.5095348,0.5267410)$ by approximately
$8.0\times10^{-5}$ in Euclidean distance. Its normalized
streamfunction extremum is $-0.1211333$, compared with
$-0.1221933$ in the reference data, a relative magnitude
difference of approximately $0.87\%$.

\begin{table}[tbp]
\centering
\caption{The vortex centers and normalized streamfunction
extrema at $\Rey=20000$ on the $2048^2$ grid, compared with the
reference streamfunction extrema of Papadopoulos
\cite{Papadopoulos2015}.}
\label{tab:trt-vortex-extrema}
\begin{tabular}{@{}lrrS[table-format=-1.7,group-digits=false]S[table-format=-1.7,group-digits=false]@{}}
\toprule
 & \multicolumn{2}{c}{Computed center} &
 \multicolumn{2}{c}{Normalized $\psi$ extremum} \\
\cmidrule(lr){2-3}\cmidrule(l){4-5}
Vortex & \multicolumn{1}{c}{$x/L$} & \multicolumn{1}{c}{$y/L$} & {PLBE-FAS} & {\cite{Papadopoulos2015}} \\
\midrule
Primary vortex  & 0.50960 & 0.52669 & -0.1211333 & -0.1221933 \\
Upper-left 1 & 0.08037 & 0.91169 &  0.0037173 &  0.0037601 \\
Upper-left 2 & 0.02440 & 0.81942 & -0.0000701 & -0.0000729 \\
Lower-left 1 & 0.04812 & 0.18277 &  0.0016282 &  0.0016399 \\
Lower-left 2 & 0.05891 & 0.05422 & -0.0000827 & -0.0000856 \\
Lower-right 1 & 0.72354 & 0.04334 &  0.0027895 &  0.0028081 \\
Lower-right 2 & 0.93031 & 0.10447 & -0.0004563 & -0.0004645 \\
\bottomrule
\end{tabular}
\end{table}

Finally, the separation and reattachment points of the resolved
corner vortices are identified from the zeros of the
wall shear stress. Their normalized coordinates are
reported in \cref{tab:cavity-reattachment}.

\begin{table}[htbp]
\centering
\caption{Major corner-vortex reattachment locations in the lid-driven cavity. Values are rounded to $0.001$.}
\label{tab:cavity-reattachment}
\begin{tabular}{lll}
\hline
Corner & $\mathrm{Re}=10000$ & $\mathrm{Re}=20000$ \\
\hline
Upper-left 1 & $(0,0.644)$ & $(0,0.607)$ \\
Lower-left 1 & $(0.352,0)$ & $(0.370,0)$ \\
Lower-right 1 & $(1,0.450)$ & $(1,0.474)$ \\
\hline
\end{tabular}
\end{table}

% Keep numerical-results floats ahead of the bibliography.
\FloatBarrier

\section{Conclusions}
\label{sec:conclusions}

In this work, we developed a nonlinear FAS multigrid solver
using the TRT-PLBE as a smoother for steady incompressible flows.
The combination of the TRT collision model and PLBE preconditioning
provides the flexibility to select relaxation parameters
with favorable smoothing properties.
Numerical experiments demonstrated improved convergence
compared with the SRT-PLBE-FAS and single-grid TRT-PLBE methods.
At high Reynolds numbers, where the V-cycle may suffer from
residual stagnation, W-cycles and hybrid V/W-cycles provide
more effective coarse-grid corrections.
Our results for 
the 2D lid-driven cavity flows at high Reynolds numbers
agree well with 
%also show good agreement with
high-resolution benchmark solutions.

The present method is developed on a hierarchy of uniform
Cartesian grids. Future work will extend the FAS framework
to locally refined Cartesian grids, with particular attention
to conservative transfer across refinement interfaces and
effective coarse-grid corrections on nonuniform mesh hierarchies.

\section*{Acknowledgements}

XZJ would like to acknowledge the support of the Science Challenge Project
under the Grant TZ2026004.

% BibTeX reads the actual cited entries from trt_plbe_fas_references.bib.
% If you add the official SIAM style file siamplain.bst to this folder,
% it is selected automatically. Until then, the bundled TeX Live `siam.bst`
% keeps this draft independently compilable (the two styles are not identical).
\IfFileExists{siamplain.bst}{%
  \bibliographystyle{siamplain}%
}{%
  \bibliographystyle{siam}%
}
\bibliography{journalnames,latticeBoltzmann,mgbib}

\end{document}